\documentclass[10pt,twoside]{article}
\usepackage{mathrsfs}
\usepackage{amsmath}
\usepackage{amssymb}
\usepackage{fancyhdr}
\usepackage{latexsym}
\usepackage{bbding}
\usepackage{mathrsfs}
\usepackage{wasysym}

\usepackage{multicol,graphics}

\newtheorem{theorem}{Theorem}[section]
\newtheorem{lemma}[theorem]{Lemma}
\newtheorem{corollary}[theorem]{Corollary}

\newtheorem{remark}[theorem]{Remark}

\numberwithin{equation}{section}
\newenvironment{proof}[1][Proof]{\noindent\textbf{#1.} }{\hfill $\Box$}
\allowdisplaybreaks

 \makeatletter
\begin{document}
\title{{A regularity criterion for the solution of the nematic liquid crystal flows in terms of $\dot{B}^{-1}_{\infty,\infty}$-norm}
\thanks{The author is partially supported by the National Natural Science Foundation of China (11171357).}
}
\author{{\small   Qiao Liu$^{1}$ \thanks{\text{E-mail address}: liuqao2005@163.com.
}}
 {\small\quad and\quad  Jihong Zhao$^{2}$ \thanks{\text{E-mail
address}: zhaojihong2007@yahoo.com.cn.}}
\\
%EndAName
{\small  $^{1}$Department of Mathematics, Hunan Normal University, Changsha, Hunan 410081,}\\
{\small People's Republic of China}\\
{\small   $^{2}$ College of Science, Northwest A\&F
University, Yangling, Shaanxi 712100,}\\
{\small    People's Republic of China}\\
}
\date{}
\maketitle

\begin{abstract}
In this paper, we investigate regularity criterion for the solution
of the nematic liquid crystal flows in dimension three and two. We
prove the solution $(u,d)$ is smooth up to time $T$ provided that
there exists a positive constant $\varepsilon_{0}>0$ such that (i)
for $n=3$,
\begin{align*}
\|(u,\nabla d)\|_{L^{\infty}(0,T;\dot{B}^{-1}_{\infty,\infty})}\leq
\varepsilon_{0},
\end{align*}
and (ii) for $n=2$,
\begin{align*}
\|\nabla d\|_{L^{\infty}(0,T;\dot{B}^{-1}_{\infty,\infty})}\leq
\varepsilon_{0}.
\end{align*}
\medskip


\end{abstract}

\section{Introduction}\label{Int}

\noindent

Liquid crystal, which is a state of matter capable of flow, but its
molecules may be oriented in a crystal-like way. Liquid crystals
exhibit a phase of matter that has properties between those of a
conventional liquid and those of a solid crystal, hence, it is
commonly considered as the fourth state of matter, different from
gases, liquid, and solid. There have been numerous attempts to
formulate continuum theories describing the behaviour of liquid
crystals flows, we refer to the seminal papers \cite{ER,LE}.  To the
present state of knowledge, three main types of liquid crystals are
distinguished, nematic, termed smectic and cholesteric. The nematic
phase appears to be the most common, where the molecules do not
exhibit any positional order, but they have long-range orientational
order.

In the present paper, we consider the following hydrodynamic model
for the flow of the nematic liquid crystal material in
$n$-dimensions ($n= 2$ or $3$):
\begin{align}
   \label{eq1.1}
&{\partial_{t}}u-\nu\Delta u +(u\cdot\nabla)u+\nabla{P}=-\lambda\nabla\cdot(\nabla d \odot\nabla d)\quad\text{ in }\mathbb{R}^{n}\times (0,+\infty),\\
%--------------------(eq1.1)------------------------------------
   \label{eq1.2}
&\partial_{t}d+(u\cdot\nabla)d=\gamma(\Delta d+|\nabla
d|^{2}d)\quad\quad
  \quad\quad\quad\quad\quad\quad\text{ in }\mathbb{R}^{n}\times (0,+\infty),\\
%--------------------(eq1.2)------------------------------------
   \label{eq1.3}
&\nabla\cdot u=0\qquad\quad\quad\quad\quad\quad\quad\quad\quad
\quad\quad\quad\qquad\quad\quad\quad\text{ in }\mathbb{R}^{n}\times (0,+\infty),\\
%--------------------(eq1.3)------------------------------------
   \label{eq1.4}
&(u,d)|_{t=0}=(u_{0},d_{0})\quad\quad\quad\quad\quad\quad\quad\quad\quad\quad\quad\quad\quad\quad
\text{in }\mathbb{R}^{n},
%--------------------(eq1.4)------------------------------------
\end{align}
where $u(x,t):\mathbb{R}^{n}\times (0,+\infty)\rightarrow
\mathbb{R}^{n}$ is the unknown velocity field of the flow,
$d:\mathbb{R}^{n}\times (0,+\infty)\rightarrow \mathbb{S}^{2}$, the
unit sphere in $\mathbb{R}^{3}$, is the unknown (averaged)
macroscopic/continuum molecule orientation of the nematic liquid
crystal flow and $P(x,t):\mathbb{R}^{n}\times (0,+\infty)\rightarrow
\mathbb{R}$ is the scalar pressure, $\nabla \cdot u=0$ represents
the incompressible condition, $u_{0}$ is a given initial velocity
with $\nabla \cdot u_{0}=0$ in distribution sense, $d_{0}:
\mathbb{R}^{n}\rightarrow \mathbb{S}^{2}$ is a given initial liquid
crystal orientation field, and the constants $\nu,\lambda,\gamma$
are positive constants that represent viscosity, the competition
between kinetic energy and potential energy,  microscopic elastic
relaxation time for the molecular orientation field. The notation
$\nabla d\odot\nabla d$ denotes the $n\times n$ matrix whose
$(i,j)$-th entry is given by $\partial_{i}d\cdot
\partial_{j}d$ ($1\leq i,j\leq n$), and there holds
$\nabla\cdot(\nabla d\odot \nabla d)= \Delta d\cdot\nabla
d+\frac{1}{2}\nabla |\nabla d|^{2}$. Since the concrete values of
the constants $\nu$, $\lambda$ and $\gamma$ do not play a special
role in our discussion, for simplicity, we assume that they all
equal to one throughout this paper.

The above system \eqref{eq1.1}--\eqref{eq1.4}  is a simplified
version of the Ericksen--Leslie model for the hydrodynamics of the
nematic liquid crystals developed during the period of 1958 through
1968 (see \cite{ER,LE}). It can be viewed as the incompressible
Navier--Stokes equations (the case $d\equiv 1$, see
\cite{BCD,YG,KT,PG}) coupling the heat flow of a harmonic map (the
case $u\equiv 0$, see \cite{CDY,CS,PG,LD,W}). The current form of
system \eqref{eq1.1}--\eqref{eq1.4} was first proposed by Lin
\cite{L} back in the late 1980's. For the mathematical analysis of
system \eqref{eq1.1}--\eqref{eq1.4}, Lin, Lin and Wang in \cite{LLW}
established that there exists global Leray--Hopf type weak solutions
to the initial boundary value problem for system
\eqref{eq1.1}--\eqref{eq1.4} on bounded domains in two space
dimensions (see also \cite{HMC}). The uniqueness of such weak
solutions is proved by Lin and Wang \cite{LW}. Wen and Ding
\cite{WD} obtained local existence and uniqueness of strong
solution.  Wang in \cite{W} proved that if the initial data
$(u_{0},d_{0})\in BMO^{-1}\times BMO$ is sufficiently small, then
system \eqref{eq1.1}--\eqref{eq1.4} has a global mild solution.
Recently, Hinmeman and Wang \cite{JHW} established the global
well--posedness of system \eqref{eq1.1}--\eqref{eq1.4} in dimension
three with small initial data $(u_{0},d_{0})$ in $L^{3}_{uloc}$,
where $L^{3}_{uloc}$ is the space of uniformly locally
$L^{3}$-integrable functions in $\mathbb{R}^{3}$. When the term
$|\nabla d|^{2}$ in \eqref{eq1.2} is replaced by
$\frac{1-|d|^{2}}{\varepsilon^{2}}$, system
\eqref{eq1.1}--\eqref{eq1.4} becomes the Ginzburg--Landau
approximation of the simplified Ericksen--Leslie system. Lin and Liu
\cite{LL1} proved the local existence of classical solutions and the
global existence of weak solutions in dimension two and three with
Dirichlet boundary conditions. For any fixed $\varepsilon$, they
also obtained the existence and uniqueness of global classical
solution either in dimension two or three for large fluid viscosity.
Later, in \cite{LL2}, they further proved that the one-dimensional
spacetime Hausdorff measure of the singular set of the so-called
suitable weak solutions is zero. For more researches about system
\eqref{eq1.1}--\eqref{eq1.4}, we refer the readers to
\cite{HW,XLW,LD,LNW,SL,WD} and the references therein.

In this paper, we are interested in the local-in-time classical
solution to system \eqref{eq1.1}--\eqref{eq1.4}. Since the strong
solutions of  the heat flow of  harmonic maps must be blowing up at
finite time \cite{CDY}, we cannot expect that
\eqref{eq1.1}--\eqref{eq1.4} has a global smooth solution with
general initial data. By using standard methods, it is known that if
the initial velocity $u_{0}\in H^{s}(\mathbb{R}^{n},\mathbb{R}^{n})$
with $\nabla\cdot u_{0}$ and $d_{0}\in
H^{s+1}(\mathbb{R}^{n},\mathbb{S}^{2})$ with $s\geq n$, then there
exists $0<T_{*}<+\infty$ depending only on the initial value such
that the system \eqref{eq1.1}--\eqref{eq1.4} has a unique local
classical solution $(u,d)$ satisfying (see for example \cite{WD})
\begin{align}\label{eq1.5}
&u\in C([0,T_{*});H^{s}(\mathbb{R}^{n},\mathbb{R}^{n}))\cap
C^{1}([0,T_{*});H^{s-1}(\mathbb{R}^{n},\mathbb{R}^{n}))\quad
\text{and
}\nonumber\\
& d\in C([0,T_{*});H^{s+1}(\mathbb{R}^{n},\mathbb{S}^{2}))\cap
C^{1}([0,T_{*});H^{s}(\mathbb{R}^{n},\mathbb{S}^{2})).
\end{align}
%--------------------(eq1.5)------------------------------------
%for all $0<T<T_{*}$.
Here, we emphasize that  such an existence theorem gives no
indication as to whether solutions actually lose their regularity or
the manner in which they may do so. %Assume that such $T_{*}$ is the
%maximum value for \eqref{eq1.5} holds, i.e., $T_{*}$ is the first
%finite singular time, the purpose of this paper is to give some
%criterion that characterizes such a $T_{*}$.
Assume that $(0,T_{*})$ is the interval for \eqref{eq1.5} holds, the
purpose of this paper is to give some criterion to ensure the
solution $(u,d)$ is smooth up to time $T_{*}$.

For the well-known Navier--Stokes equations with dimension $n\geq
3$, there are many interesting sufficient conditions for regularity
of solutions (see for example \cite{HBV,GG11,KOT,KT}), and the
Ladyzhenskaya--Prodi--Serrin condition (see \cite{JS,PG}) state
that %the condition
\begin{align}\label{eq1.6}
u\in L^{\alpha}(0,T_{*};L^{\beta}(\mathbb{R}^{n})) \text{ for all }
\frac{2}{\alpha}+\frac{n}{\beta}\leq 1, 2\leq \alpha<\infty,
n<\beta\leq\infty
\end{align}
%--------------------(eq1.6)------------------------------------
ensure the smoothness of solution $u$ up to time $T_{*}$. The
limiting case $u\in L^{\infty}(0,T_{*};L^{n}(\mathbb{R}^{n}))$ in
\eqref{eq1.6} has been proved by Escauriaza, Seregin and
\u{S}ver\'{a}k \cite{ESS} by using the method of backward uniqueness
of solution.  Beale, Kato and Majda in \cite{BKM} proved that the
vorticity $\omega=\nabla\times u$ does not belong to
$L^{1}(0,T_{*};L^{\infty}(\mathbb{R}^{n}))$ if $T_{*}$ is the first
finite singular time. On the other hand, as for the heat flow of
harmonic maps into $\mathbb{S}^{2}$, Wang \cite{W08} established
that for $n\geq 2$, the condition $\nabla d\in
L^{\infty}(0,T;L^{n}(\mathbb{R}^{n}))$ implies that the solution $d$
is regular on $(0,T]$, i.e., $d\in C^{\infty}((0,T]\times
\mathbb{R}^{n})$. When $n=2$, Lin, Lin and Wang obtained that the
local smooth solution $(u,d)$ to \eqref{eq1.1}--\eqref{eq1.4} can be
continued past any time $T>0$ provided that there holds
\begin{align*}
\int_{0}^{T}\|\nabla d (\cdot, t)\|_{L^{4}}^{4}\text{d}t<\infty.
\end{align*}
Huang and Wang \cite{HW1} established that if $0<T_{*}<\infty$ is
the first finite singular time of the smooth solutions $(u,d)$ to
system \eqref{eq1.1}--\eqref{eq1.4}, then
\begin{align*}
&\int_{0}^{T_{*}} (\|\omega\|_{L^{\infty}}+\|\nabla
d\|_{L^{\infty}}^{2})\text{d}x=\infty \quad \text{ when dimension }
n=3;\\
&\quad\quad\int_{0}^{T_{*}} \|\nabla
d\|_{L^{\infty}}^{2}\text{d}x=\infty \quad \text{ when dimension }
n=2.
\end{align*}
In the references cited above, we noticing that the scaling
invariance property plays a particularly significant role. For
system \eqref{eq1.1}--\eqref{eq1.4}, it is clear that if
$(u(x,t),d(x,t))$ is the solution of system
\eqref{eq1.1}--\eqref{eq1.4}, then
\begin{align}\label{eq1.7}
(u_{\lambda}(x,t),d_{\lambda}(x,t)):= (\lambda u (\lambda
x,\lambda^{2}t),d (\lambda x,\lambda^{2}t))
\end{align}
%--------------------(eq1.7)------------------------------------
for any $\lambda>0$ is also the solution of
\eqref{eq1.1}--\eqref{eq1.3} with initial data
$(u_{0\lambda}(x),d_{0\lambda}(x)):=(\lambda u_{0}(\lambda
x),d_{0}(\lambda x))$. %From the scaling invariant point of view,
%when space dimension $n=2$, it is easy to see that the space
%$L^{4}(0,T;L^{4}(\mathbb{R}^{2}))$ is scaling invariant for
%$(u,\nabla d)$, that is
%\begin{align*}
%\|(u_{\lambda},\nabla d_{\lambda})\|_{L^{4}(0,T;L^{4})}=\|(u,\nabla
%d)\|_{L^{4}(0,T;L^{4})}.
%\end{align*}
%Similarly,  when $n=3$, the space
%$L^{1}(0,T;L^{\infty}(\mathbb{R}^{3}))$ is scaling invariant for
%$(\omega,\nabla d)$, that is
%\begin{align*}
%\|(\omega_{\lambda},\nabla
%d_{\lambda})\|_{L^{1}(0,T;L^{\infty})}=\|(\omega,\nabla
%d)\|_{L^{1}(0,T;L^{\infty})}.
%\end{align*}
%Here $\omega_{\lambda}:=\omega_{\lambda}(x,t)=\lambda^{2}
%\omega(\lambda x,\lambda^{2}t)$.
In fact,  it is easy to verify that the space
$L^{n}(\mathbb{R}^{n})\times
\dot{W}^{1,n}(\mathbb{R}^{n})$\footnote{Here
$\dot{W}^{1,n}(\mathbb{R}^{n})$ denotes the homogeneous Sobolev
space on $\mathbb{R}^{n}$ (see e.g., \cite{LD}).} is the scaling
invariance space for $(u(t),d(t))$ in system
\eqref{eq1.1}--\eqref{eq1.4}, i.e., $L^{n}(\mathbb{R}^{n})\times
L^{n}(\mathbb{R}^{n})$-norm of $(u(t),\nabla d(t))$ is invariant
under the action of the scaling \eqref{eq1.7}.  Due to the facts
that
\begin{align*}
L^{n}(\mathbb{R}^{n})\subset
\dot{B}^{-1}_{\infty,\infty}(\mathbb{R}^{n}) \text{ and }
L^{n}(\mathbb{R}^{n})\neq
\dot{B}^{-1}_{\infty,\infty}(\mathbb{R}^{n}),
\end{align*}
and from a mathematical viewpoint, Besov space
$\dot{B}^{-1}_{\infty,\infty}(\mathbb{R}^{n})$ is the largest
scaling invariant space of the system \eqref{eq1.1}--\eqref{eq1.4},
% which include general Sobolev spaces,
the purpose of this paper is to establish a regularity criterion for
local-in-time smooth solutions of system
\eqref{eq1.1}--\eqref{eq1.4} in term of the homogeneous Besov space
$\dot{B}^{-1}_{\infty,\infty}$-norm.

Our main results are as follows:

\begin{theorem}\label{thm1.1}
For $n=3$, $u_{0}\in H^{3}(\mathbb{R}^{3},\mathbb{R}^{3})$ with
$\nabla\cdot u_{0}=0$ and $d_{0}\in
H^{4}(\mathbb{R}^{3},\mathbb{S}^{2})$. Let $0<T_{*}<\infty$ be the
 value such that the nematic liquid crystal flow
\eqref{eq1.1}--\eqref{eq1.4} has a unique solution $(u,d)$
satisfying \eqref{eq1.5}. If there exists a small positive constant
$\varepsilon_{0}$ such that
\begin{align}\label{eq1.8}
\|(u,\nabla
d)\|_{L^{\infty}(0,T_{*};\dot{B}^{-1}_{\infty,\infty})}\leq
\varepsilon_{0},
\end{align}
%--------------------(eq1.8)------------------------------------
then $(u,d)$ is smooth up to time $t=T_{*}$.
\end{theorem}
%--------------------(thm1.1)------------------------------------

\begin{remark}\label{rem1.2}
In \cite{ESS}, Escauriaza, Seregin and \u{S}ver\'{a}k used the fact
that functions in $L^{3}(\mathbb{R}^{3})$ has decay at infinity,
which ensures that the solution of the 3D Navier--Stokes equations
is smooth outside an big ball centered at origin so that the
backward uniqueness theorem can be applied. We can not generalize
the regularity criterion \eqref{eq1.8} as $(u,\nabla d)\in
L^{\infty}(0,T_{*};\dot{B}^{-1}_{\infty,\infty}(\mathbb{R}^{3}))$,
since functions in $\dot{B}^{-1}_{\infty,\infty}(\mathbb{R}^{3})$ is
different from the functions in $L^{3}(\mathbb{R}^{3})$, which has
no decay at infinity.
\end{remark}
%--------------------(rem1.2)------------------------------------

As a byproduct of our proof of Theorem 1.1, we obtain the following
corresponding criterion in dimension two. More precisely, we have

\begin{corollary}\label{cor1.2}
For n=2, $u_{0}\in H^{2}(\mathbb{R}^{2},\mathbb{R}^{2})$ with
$\nabla\cdot u_{0}=0$ and $d_{0}\in
H^{3}(\mathbb{R}^{2},\mathbb{S}^{2})$. Let $0<T_{*}<\infty$ be the
 value such that the nematic liquid crystal flow
\eqref{eq1.1}--\eqref{eq1.4} has a unique solution $(u,d)$
satisfying \eqref{eq1.5}. If there exists a small positive constant
$\varepsilon_{0}$ such that
\begin{align}\label{eq1.9}
\|\nabla d\|_{L^{\infty}(0,T_{*};\dot{B}^{-1}_{\infty,\infty})}\leq
\varepsilon_{0},
\end{align}
%--------------------(eq1.9)------------------------------------
then $(u,d)$ is smooth up to time $t=T_{*}$.
\end{corollary}
%--------------------(cor1.2)------------------------------------

The remaining parts of the paper is written as follows. Section 2,
we recall the definition of Besov spaces and an useful inequality.
Section 3 is devoted to proving  Theorem \ref{thm1.1} and Corollary
\ref{cor1.2}. Throughout the paper, $C$ denotes the positive
constant and its value may change from line to line; $\|\cdot\|_{X}$
denotes the norm of space $X(\mathbb{R}^{3})$ or
$X(\mathbb{R}^{2})$.

\section{Preliminaries and a key lemma}

In this section, we will give the definition of the Besov spaces and
an useful inequality. In order to define Besov spaces, we first
introduce the Littlewood--Paley decomposition theory. Let
$\mathcal{S}(\mathbb{R}^{n})$ be the Schwartz class of rapidly
decreasing functions, for given $f\in \mathcal{S}(\mathbb{R}^{n})$,
its Fourier transform $\mathcal{F}f=\widehat{f}$ and its inverse
Fourier transform $\mathcal{F}^{-1}f=\check{f}$ are, respectively,
defined by
\begin{align*}
\widehat{f}(\xi):=\int_{\mathbb{R}^{n}}e^{-ix\cdot\xi}f(x)\text{d}x,
\end{align*}
and
\begin{align*}
\check{f}(x)=\frac{1}{(2\pi)^{n}}\int_{\mathbb{R}^{n}}e^{ix\cdot\xi}f(x)\text{d}\xi.
\end{align*}
More generally, the Fourier transform of any given $f\in
\mathcal{S}'(\mathbb{R}^{n})$, the space of tempered distributions,
is given by
\begin{align*}
<\widehat{f},g>=<f,\widehat{g}>,\quad \text{ for any } g\in
\mathcal{S}(\mathbb{R}^{n}).
\end{align*}
%The Fourier transform is a bounded linear bijection from
%$\mathcal{S}'$ to $\mathcal{S}'$.
Let
\begin{align*}
\mathcal{S}_{h}:= \{\phi \in \mathcal{S}(\mathbb{R}^{n}),
\int_{\mathbb{R}^{n}}\phi(x)x^{\gamma}\text{d}x=0,
|\gamma|=0,1,2,\cdots \}.
\end{align*}
Then its dual is given by
\begin{align*}
\mathcal{S}'_{h}=\mathcal{S}'/\mathcal{S}^{\bot}_{h}=\mathcal{S}'/\mathcal{P},
\end{align*}
where $\mathcal{P}$ is the space of polynomial.  Let us choose two
nonnegative radial functions $\chi, \varphi \in
\mathcal{S}(\mathbb{R}^{n})$ supported in $\mathfrak{B}=\{\xi\in
\mathbb{R}^{n}:|\xi|\leq \frac{4}{3}\}$ and
$\mathfrak{C}=\{\xi\in\mathbb{R}^{n}:\frac{3}{4}\leq |\xi|\leq
\frac{8}{3}\}$ respectively, such that
\begin{align*}
\sum_{j\in\mathbb{Z}}\varphi(2^{-j}\xi)=1 \quad\text{ for any
}\xi\in\mathbb{R}^{n}\backslash\{0\},
\end{align*}
and
\begin{align*}
\chi(\xi)+\sum_{j\geq 0}\varphi(2^{-j}\xi)=1\quad\text{ for any }
\xi\in\mathbb{R}^{n}.
\end{align*}
For $j\in \mathbb{Z}$, the homogeneous Littlewood--Paley projection
operators $\dot{S}_{j}$ and $\dot{\Delta}_{j}$ are, respectively,
defined as
\begin{align*}
\dot{S}_{j}f=\chi(2^{-j}D)f=2^{nj}\int_{\mathbb{R}^{n}}\widetilde{h}(2^{j}y)f(x-y)\text{d}y,
\text{ where } \widetilde{h}=\mathcal{F}^{-1}\chi,
\end{align*}
and
\begin{align*}
 \dot{\Delta}_{j} f=\varphi(2^{-j}D)
f=2^{nj}\int_{\mathbb{R}^{n}}h(2^{-j}y)f(x-y)\text{d}y, \text{ where
} h=\mathcal{F}^{-1}\varphi.
\end{align*}
Informally, $\dot{\Delta}_{j}$ is a  frequency projection to the
annulus $\{|\xi|\sim 2^{j}\}$, while $\dot{S}_{j}$ is a frequency
projection to the ball $\{|\xi|\lesssim 2^{j} \}$. One can easily
verify that $\dot{\Delta}_{j}\dot{\Delta}_{k}f=0$ if $|j-k|\geq 2$.

Let $s\in\mathbb{R}, p,q\in [1,+\infty]$, the homogenous Besov space
$\dot{B}^{s}_{p,q}(\mathbb{R}^{n})$ is defined by those
distributions $f$ in $\mathcal{S}'_{h}$ such that %the full-dyadic
%decomposition, i.e., we say that
%$f\in\dot{B}^{s}_{p,q}(\mathbb{R}^{n})$, if $f\in \mathcal{S}'_{h}$
%and
\begin{align*}
\sum_{j\in\mathbb{Z}}(2^{js}\|\dot{\Delta}_{j} f\|_{L^{p}})^{q}<
\infty,
\end{align*}
with the norm
\begin{equation*}
\|f\|_{\dot{B}^{s}_{p,q}}:=
\begin{cases}
(\sum_{j\in\mathbb{Z}}
2^{jsq}\|\dot{\Delta}_{j}f\|_{L^{p}}^{q})^{\frac{1}{q}}, 1\leq q<+\infty,\\
\sup_{j\in\mathbb{Z}} \{2^{js}\|\dot{\Delta}_{j}f\|_{L^{p}}\},
q=+\infty.
\end{cases}
\end{equation*}

The following interpolation inequality will be used in the next
section.

\begin{lemma}\label{lem2.1}
(see \cite{BCD,CPG}) Let $1\leq q<p<\infty$ and $\alpha$ be a
positive real number. A constant $C$ exists such that
\begin{align*}
\|f\|_{L^{p}}\leq
C\|f\|_{\dot{B}^{-\alpha}_{\infty,\infty}}^{1-\theta}\|f\|_{\dot{B}^{\beta}_{q,q}}^{\theta}
\text{ with }\beta=\alpha(\frac{p}{q}-1) \text{ and }
\theta=\frac{q}{p},
\end{align*}
for all $f\in \dot{B}^{-\alpha}_{\infty,\infty}(\mathbb{R}^{n})\cap
\dot{B}^{\beta}_{q,q}(\mathbb{R}^{n}) $ with $n\geq 1$.
\end{lemma}
%--------------------(lem2.1)------------------------------------

It is of interest to notice that the homogeneous Besov space
$\dot{B}^{s}_{2,2}(\mathbb{R}^{n})$ is equivalent to the homogeneous
Sobolev space $\dot{H}^{s}(\mathbb{R}^{n})$. Hence, from Lemma
\ref{lem2.1} above, we have the following interpolation inequality:
\begin{align}\label{eq2.1}
\|f\|_{L^{q}}\leq
C\|f\|_{\dot{H}^{\alpha(\frac{q}{2}-1)}}^{\frac{2}{q}}\|f\|_{\dot{B}^{-\alpha}_{\infty,\infty}}^{1-\frac{2}{q}}
\text{ with }2<q<\infty \text{ and } \alpha>0,
\end{align}
%--------------------(eq2.1)------------------------------------
for all $f\in \dot{H}^{\alpha(\frac{q}{2}-1)}(\mathbb{R}^{n})\cap
\dot{B}^{-\alpha}_{\infty,\infty}(\mathbb{R}^{n}) $ with $n\geq 1$.

\section{The proofs of Theorem \ref{thm1.1} and Corollary \ref{cor1.2}}

In this section, we shall give the proofs of Theorem \ref{thm1.1}
and Corollary \ref{cor1.2}.  We first need to prove the following
lemma.

\begin{lemma}\label{lem3.1}
For $n=2$ or $3$, $s\geq n$, $u_{0}\in
H^{s}(\mathbb{R}^{n},\mathbb{R}^{n})$ with $\nabla \cdot u_{0}=0$
and $d_{0}\in H^{s+1}(\mathbb{R}^{n},\mathbb{S}^{2})$,
$0<T_{*}<\infty$, let $(u,d)$ be a solution to system
\eqref{eq1.1}--\eqref{eq1.4} satisfying \eqref{eq1.5}, and there
exists a small positive constant $\varepsilon_{0}$ such that
\begin{align}\label{eq3.1}
&\|(u,\nabla
d)\|_{L^{\infty}(0,T_{*};\dot{B}^{-1}_{\infty,\infty})}\leq
\varepsilon_{0},\quad\quad \text{ for } n=3;
\end{align}
%--------------------(eq3.1)------------------------------------
or
\begin{align}\label{eq3.2}
&\|\nabla d\|_{L^{\infty}(0,T_{*};\dot{B}^{-1}_{\infty,\infty})}\leq
\varepsilon_{0}, \quad\quad \text{ for } n=2.
\end{align}
%--------------------(eq3.2)------------------------------------
Then
\begin{align}\label{eq3.3}
\sup_{0<t\leq T_{*}}\|\nabla u(\cdot,t)\|_{L^{2}}^{2}+\|\Delta
d(\cdot,t)\|_{L^{2}}^{2}+\int_{0}^{T_{*}}\left(\|\nabla^{2}u(\cdot,t)\|_{L^{2}}^{2}+\|\nabla\Delta
d(\cdot,t)\|_{L^{2}}^{2}\right)\text{d}t \leq C_{0},
\end{align}
%--------------------(eq3.3)------------------------------------
where $C_{0}$ is a positive constant depending only on $u_{0}$,
$d_{0}$ and $T_{*}$.
\end{lemma}
%--------------------(lem3.1)------------------------------------

\begin{proof}\label{proof of lem3.1}
We firstly notice that for all smooth solutions to system
\eqref{eq1.1}--\eqref{eq1.4}, one has the following basic energy law
(see \cite{LLW,LL1}):
\begin{align}\label{eq3.4}
\|u(t)\|_{L^{2}}^{2}+&\|\nabla d
(t)\|_{L^{2}}^{2}+\int_{0}^{t}(\|\nabla
u(\tau)\|_{L^{2}}^{2}+\|(\Delta
d+|\nabla d|^{2} )(\tau)\|_{L^{2}}^{2})\text{d}\tau\nonumber\\
\leq& \|u_{0}\|_{L^{2}}^{2}+\|\nabla d_{0}\|_{L^{2}}^{2},\quad\text{
for all } 0<t<\infty,
\end{align}
%--------------------(eq3.4)------------------------------------
and when the space dimension $n=2$, the above energy inequality
becomes energy equality.

Now, applying $\nabla$ to the equation \eqref{eq1.1}, multiplying
the resulting equation by $\nabla u$, integrating with respect to
$x$ over $\mathbb{R}^{n}$ with $n=2$ or $3$, and using integration
by parts, we get
\begin{align}\label{eq3.5}
\frac{1}{2}\frac{d}{dt}\|\nabla u(t)\|_{L^{2}}^{2}+\|\Delta
u(t)\|_{L^{2}}^{2}=-\int_{\mathbb{R}^{n}} \nabla (u\cdot\nabla
u)\nabla u\text{d}x-\int_{\mathbb{R}^{n}}\nabla(\nabla\cdot(\nabla
d\odot\nabla d))\nabla u\text{d}x.
\end{align}
%--------------------(eq3.5)------------------------------------
Similarly, applying $\nabla^{2}$ to the equation \eqref{eq1.2},
multiplying the resulting equation by $\nabla^{2} d$, integrating
with respect to $x$ over $\mathbb{R}^{n}$ with $n=2$ or $3$, and
using integration by parts, we get
\begin{align}\label{eq3.6}
\frac{1}{2}\frac{d}{dt}\|\nabla^{2} d(t)\|_{L^{2}}^{2}+\|\nabla
\Delta
d(t)\|_{L^{2}}^{2}=-\int_{\mathbb{R}^{n}}\nabla^{2}(u\cdot\nabla
d)\nabla^{2}d\text{d}x +\int_{\mathbb{R}^{n}}\nabla^{2}(|\nabla
d|^{2}d)\nabla^{2}d\text{d}x.
\end{align}
%--------------------(eq3.6)------------------------------------
Combining \eqref{eq3.5} and \eqref{eq3.6} together, and using the
fact $\nabla\cdot u=0$, we get
\begin{align}\label{eq3.7}
&\frac{1}{2}\frac{d}{dt}(\|\nabla u(t)\|_{L^{2}}^{2}+\|\nabla^{2}
d(t)\|_{L^{2}}^{2})+(\|\Delta u(t)\|_{L^{2}}^{2}+\|\nabla \Delta
d(t)\|_{L^{2}}^{2})\nonumber\\
=&-\int_{\mathbb{R}^{n}}[\nabla(u\cdot\nabla u)-u\cdot\nabla\nabla
u]\nabla u\text{d}x-\int_{\mathbb{R}^{n}}\nabla(\nabla\cdot(\nabla
d\odot\nabla d))\nabla u\text{d}x\nonumber\\
&-\int_{\mathbb{R}^{n}}[\nabla^{2}(u\cdot\nabla
d)-u\cdot\nabla\nabla^{2}d]\nabla^{2}d\text{d}x+\int_{\mathbb{R}^{n}}\nabla^{2}(|\nabla
d|^{2} d)\nabla^{2} d\text{d}x\nonumber\\
=&\!-\!\!\int_{\!\mathbb{R}^{n}}\![\nabla u\cdot\nabla u\nabla
u\!+\!\nabla(\nabla\cdot(\nabla d\odot\nabla d))\nabla
u\!+\!\nabla^{2} u\cdot\nabla d\nabla^{2}d\!+\!2\nabla
u\cdot\nabla\nabla d\nabla^{2}d\!-\!\nabla^{2}(|\nabla
d|^{2}d)\nabla^{2}d]\text{d}x\nonumber\\
\triangleq& I_{1}+I_{2}+I_{3}+I_{4}+I_{5}.
\end{align}
%--------------------(eq3.7)------------------------------------
By using the H\"{o}lder's inequality,  the interpolation inequality
and Lemma \ref{lem2.1}, we obtain
\begin{align}\label{eq3.8}
I_{1}=&\int_{\mathbb{R}^{n}}u\cdot\nabla u\nabla^{2}u\text{d}x\nonumber\\
\leq&
\begin{cases}
C\|u\|_{L^{6}}\|\nabla u\|_{L^{3}}\|\nabla^{2} u\|_{L^{2}}  \leq
C\|u\|_{\dot{H}^{2}}^{\frac{1}{3}}\|u\|_{\dot{B}^{-1}_{\infty,\infty}}^{\frac{2}{3}}
(\|\nabla u\|_{\dot{H}^{1}}^{\frac{2}{3}}\|\nabla
u\|_{\dot{B}^{-2}_{\infty,\infty}}^{\frac{1}{3}})
\|\Delta u\|_{L^{2}}\\
\leq C\|u\|_{\dot{B}^{-1}_{\infty,\infty}}\|\Delta u\|_{L^{2}}^{2}, \quad\quad\text{ for } n=3\\
C\|u\|_{L^{4}}\|\nabla u\|_{L^{4}}\|\nabla^{2} u\|_{L^{2}}\leq C
\|u\|_{L^{2}}^{\frac{1}{2}}\|\nabla
u\|_{L^{2}}^{\frac{1}{2}}(\|\nabla u\|_{L^{2}}^{\frac{1}{2}}\|\Delta
u\|_{L^{2}}^{\frac{1}{2}})\|\Delta u\|_{L^{2}}\\
\leq \frac{1}{16}\|\Delta u\|_{L^{2}}^{2}+C\|u\|_{L^{2}}\|\nabla
u\|_{L^{2}}^{2} \quad (\text{ by using energy inequality \eqref{eq3.4}} )\\
\leq \frac{1}{16}\|\Delta u\|_{L^{2}}^{2}+C \|\nabla
u\|_{L^{2}}^{2}, \quad\quad \text{ for } n=2.
\end{cases}
\end{align}
%--------------------(eq3.8)------------------------------------
\begin{align}\label{eq3.9}
I_{2}=&\int_{\mathbb{R}^{n}}\nabla\cdot(\nabla d\odot\nabla
d)\nabla^{2}u\text{d}x=\int_{\mathbb{R}^{n}} [\nabla d\Delta
d+\nabla(\frac{|\nabla
d|^{2}}{2})]\nabla^{2}u\text{d}x=\int_{\mathbb{R}^{n}} \nabla
d\Delta d \nabla^{2}u\text{d}x\nonumber\\
\leq&
\begin{cases}
C\|\nabla d\|_{L^{6}}\|\Delta d\|_{L^{3}}\|\nabla^{2} u\|_{L^{2}}
\leq C \|\nabla d\|_{\dot{H}^{2}}^{\frac{1}{3}}\|\nabla
d\|_{\dot{B}^{-1}_{\infty,\infty}}^{\frac{2}{3}}(\|\Delta
d\|_{\dot{H}^{1}}^{\frac{2}{3}}\|\Delta
d\|_{\dot{B}^{-2}_{\infty,\infty}}^{\frac{1}{3}})\|\Delta
u\|_{L^{2}}\\
\leq C\|\nabla d\|_{\dot{B}^{-1}_{\infty,\infty}}\|\nabla \Delta
d\|_{L^{2}}\|\Delta
u\|_{L^{2}}\\
\leq C\|\nabla d\|_{\dot{B}^{-1}_{\infty,\infty}}(\|\nabla \Delta
d\|_{L^{2}}^{2}+\|\Delta
u\|_{L^{2}}^{2}),\quad \quad \text{ for } n=3;\\
C\|\nabla d\|_{L^{4}}\|\Delta d\|_{L^{4}}\|\nabla^{2}
u\|_{L^{2}}\leq C\|\nabla d\|_{L^{2}}^{\frac{1}{2}}\|\nabla^{2}
d\|_{L^{2}}^{\frac{1}{2}}(\|\nabla^{2}
d\|_{L^{2}}^{\frac{1}{2}}\|\nabla\Delta
d\|_{L^{2}}^{\frac{1}{2}})\|\Delta u\|_{L^{2}}\\
\leq \frac{1}{16}(\|\Delta u\|_{L^{2}}^{2}+\|\nabla \Delta
d\|_{L^{2}}^{2})+C\|\nabla d\|_{L^{2}}\|\nabla^{2}
d\|_{L^{2}}^{2}\quad (\text{ by \eqref{eq3.4}} )\\
\leq \frac{1}{16}(\|\Delta u\|_{L^{2}}^{2}+\|\nabla \Delta
d\|_{L^{2}}^{2})+C\|\nabla^{2} d\|_{L^{2}}^{2}, \quad\quad \text{
for } n=2.
\end{cases}
\end{align}
%--------------------(eq3.9)------------------------------------
Similar as the estimates of $I_{2}$, we obtain
\begin{align}\label{eq3.10}
\!I_{3}%=-\int_{\mathbb{R}^{n}} \nabla^{2} u \cdot \nabla d\nabla^{2}
%d\text{d}x
\!\!\leq\!
\begin{cases}
C\|\nabla d\|_{L^{6}}\|\Delta d\|_{L^{3}}\|\nabla^{2}
u\|_{L^{2}}\leq C\|\nabla
d\|_{\dot{B}^{-1}_{\infty,\infty}}(\|\nabla \Delta
d\|_{L^{2}}^{2}+\|\Delta
u\|_{L^{2}}^{2}),\quad \text{ for } n=3;\\
C\|\nabla d\|_{L^{4}}\|\Delta d\|_{L^{4}}\|\nabla^{2}
u\|_{L^{2}}\leq  \frac{1}{16}(\|\Delta u\|_{L^{2}}^{2}+\|\nabla
\Delta d\|_{L^{2}}^{2})+C\|\nabla^{2} d\|_{L^{2}}^{2}, \quad \text{
for } n=2.
\end{cases}
\end{align}
%--------------------(eq3.10)------------------------------------
For the term $I_{4}$, we have
\begin{align}\label{eq3.11}
I_{4}\!\leq\!
\begin{cases}
C\|\nabla u\|_{L^{3}}\|\nabla^{2} d\|_{L^{3}}^{2}\leq
C\|u\|_{L^{2}}^{\frac{1}{4}}\|\Delta
u\|_{L^{2}}^{\frac{3}{4}}\|\nabla^{2}
d\|_{\dot{H}^{1}}^{\frac{4}{3}}\|\nabla^{2}
d\|_{\dot{B}^{-2}_{\infty,\infty}}^{\frac{2}{3}}\\
\leq\! C\|\Delta u\|_{L^{2}}^{\frac{3}{4}}\|\nabla\Delta
d\|_{L^{2}}^{\frac{4}{3}}\|\nabla
d\|_{\!\dot{B}^{-1}_{\!\infty,\!\infty}}^{\frac{2}{3}}\! \quad \text{ by \eqref{eq3.4}}\\
\leq \frac{1}{4}\|\Delta u\|_{L^{2}}^{2}\!+C\|\nabla
d\|_{\!\dot{B}^{-1}_{\!\infty,\!\infty}}
\|\nabla\Delta d\|_{L^{2}}^{2}+C, \text{ for } n=3;\\
C\|\nabla u\|_{L^{3}}\|\nabla^{2} d\|_{L^{3}}^{2}\leq  C\|
u\|_{L^{2}}^{\frac{1}{3}}\|\Delta u\|_{L^{2}}^{\frac{2}{3}}
(\|\nabla^{2}d\|_{\dot{H}^{1}}^{\frac{2}{3}}\|\nabla^{2}
d\|_{\dot{B}^{-2}_{\infty,\infty}}^{\frac{1}{3}})^{2}\\
\leq C \|\Delta u\|_{L^{2}}^{\frac{2}{3}}\|\nabla\Delta
d\|_{L^{2}}^{\frac{4}{3}}
\|\nabla d\|_{\dot{B}^{-1}_{\infty,\infty}}^{\frac{2}{3}} \quad \text{ by \eqref{eq3.4}}\\
\leq \frac{1}{16}\|\Delta u\|_{L^{2}}^{2}+C\|\nabla
d\|_{\dot{B}^{-1}_{\infty,\infty}}\|\nabla\Delta d\|_{L^{2}}^{2} +C,
\quad \text{ for } n=2,
\end{cases}
\end{align}
%--------------------(eq3.11)------------------------------------
where we have used the following Gagliardo--Nirenberg inequality
\begin{align*}
\|\nabla u\|_{L^{3}}\leq C\|u\|_{L^{2}}^{\frac{1}{4}}\|\Delta
u\|_{L^{2}}^{\frac{3}{4}}, \text{ when } n=3;\\
\|\nabla u\|_{L^{3}}\leq C\|u\|_{L^{2}}^{\frac{1}{3}}\|\Delta
u\|_{L^{2}}^{\frac{2}{3}}, \text{ when } n=2.
\end{align*}
By using  facts $|d|=1$ and $\Delta d\cdot d=-|\nabla d|^{2} $, we
see that
\begin{align}\label{eq3.12}
I_{5}&=\int_{\mathbb{R}^{n}}\nabla^{2}(|\nabla d|^{2}
d)\nabla^{2}d\text{d}x\nonumber\\
&=\int_{\mathbb{R}^{n}}\nabla(2\nabla^{2}d d+|\nabla d|^{2}
d)\nabla^{2}d\text{d}x\nonumber\\
&=\int_{\mathbb{R}^{n}}(2\nabla^{3}d\nabla d
d+2|\nabla^{2}d|^{2}d+5\nabla^{2} d|\nabla
d|^{2})\nabla^{2}d\text{d}x\nonumber\\
&=\int_{\mathbb{R}^{n}}(2\nabla^{3}d\nabla d
d+2|\nabla^{2}d|^{2}d+5\nabla^{2} d\Delta d d)\nabla^{2}d\text{d}x\nonumber\\
 &\leq C\int_{\mathbb{R}^{n}}(|\nabla^{3} d||\nabla^{2}d ||\nabla
d|+|\nabla^{2} d|^{3}+|\nabla^{2}d||\Delta d|)\text{d}x\nonumber\\
&\leq C\|\nabla d\|_{L^{6}}\|\nabla^{2} d\|_{L^{3}}\|\nabla\Delta
d\|_{L^{2}}
+C\|\nabla^{2} d\|_{L^{3}}^{3}+C\|\nabla^{2} d\|_{L^{3}}^{2}\|\Delta d\|_{L^{3}}\nonumber\\
&\leq\! C\|\nabla\Delta d\|_{\!L^{\!2}}(\|\nabla
d\|_{\!\dot{H}^{\!2}}^{\frac{1}{3}} \|\nabla
d\|_{\!\dot{B}^{-1}_{\!\infty,\!\infty}}^{\frac{2}{3}}\!)
(\|\nabla^{2} d\|_{\!\dot{H}^{1}}^{\frac{2}{3}}\|\nabla^{2} d
\|_{\!\dot{B}^{-2}_{\!\infty,\!\infty}}^{\frac{1}{3}}\!)\!
+\!C\!(\|\nabla^{2} d\|_{\!\dot{H}^{1}}^{\frac{2}{3}}
\|\nabla^{2} d\|_{\!\dot{B}^{-2}_{\!\infty,\!\infty}}^{\frac{1}{3}})^{3}\nonumber\\
&\leq C\|\nabla d\|_{\dot{B}^{-1}_{\infty,\infty}}\|\nabla\Delta
d\|_{L^{2}}^{2},\quad\quad \text{ for } n=2 \text{ or }  3.
\end{align}
%--------------------(eq3.12)------------------------------------
Inserting \eqref{eq3.8}--\eqref{eq3.12} into \eqref{eq3.7}, one gets
\begin{align*}
&\frac{1}{2}\frac{d}{dt}(\|\nabla u(t)\|_{L^{2}}^{2}+\|\nabla^{2}
d(t)\|_{L^{2}}^{2})+(\|\Delta u(t)\|_{L^{2}}^{2}+\|\nabla \Delta
d(t)\|_{L^{2}}^{2})\nonumber\\
\leq&
\begin{cases}
\frac{1}{2}\|\Delta u\|_{L^{2}}^{2}\!+C(\|
u\|_{\dot{B}^{-1}_{\infty,\infty}}+\|\nabla
d\|_{\!\dot{B}^{-1}_{\!\infty,\!\infty}})
(\|\Delta u\|_{L^{2}}^{2}+\|\nabla\Delta d\|_{L^{2}}^{2})+C\\
=\frac{1}{4}\|\Delta u\|_{L^{2}}^{2}\!+C\|( u,\nabla
d)\|_{\dot{B}^{-1}_{\infty,\infty}}(\|\Delta
u\|_{L^{2}}^{2}+\|\nabla\Delta d\|_{L^{2}}^{2})+C\\
 \leq
\frac{1}{4}\|\Delta u\|_{L^{2}}^{2}\!+C\varepsilon_{0}
(\|\Delta u\|_{L^{2}}^{2}+\|\nabla\Delta d\|_{L^{2}}^{2})+C,\quad \text{ for } n=3;\\
\frac{1}{4}(\|\Delta u\|_{L^{2}}^{2}+\|\nabla \Delta
d\|_{L^{2}}^{2})+C\|\nabla
d\|_{\dot{B}^{-1}_{\infty,\infty}}\|\nabla\Delta
d\|_{L^{2}}^{2}+C(\|\nabla u\|_{L^{2}}^{2}+\|\nabla^{2}
d\|_{L^{2}}^{2})+C\\
 =\frac{1}{4}(\|\Delta u\|_{L^{2}}^{2}+\|\nabla \Delta
d\|_{L^{2}}^{2})+C\varepsilon_{0}\|\nabla\Delta
d\|_{L^{2}}^{2}+C(\|\nabla u\|_{L^{2}}^{2}+\|\nabla^{2}
d\|_{L^{2}}^{2})+C, \quad \text{ for } n=2.
\end{cases}
\end{align*}
By taking the $\varepsilon_{0}$ in \eqref{eq3.1} or \eqref{eq3.2}
small enough, and noticing that there holds equality $\|\nabla^{2}
d\|_{L^{2}}^{2}=\|\Delta d\|_{L^{2}}^{2}$,  we get
\begin{align}\label{eq3.13}
&\frac{1}{2}\frac{d}{dt}(\|\nabla u(t)\|_{L^{2}}^{2}+\|\Delta
d(t)\|_{L^{2}}^{2})+\frac{1}{2}(\|\Delta u(t)\|_{L^{2}}^{2}+\|\nabla
\Delta d(t)\|_{L^{2}}^{2})\nonumber\\
\leq&
\begin{cases}
C,\quad\quad\text{ for } n=3;\\
C(\|\nabla u\|_{L^{2}}^{2}+\|\Delta d\|_{L^{2}}^{2})+C,\quad\text{
for }n=2.
\end{cases}
\end{align}
%--------------------(eq3.13)------------------------------------
Then, by integrating with respect to $t$ over $[0; T_{*}]$ for
$n=3$, or by using the Gronwall's inequality for $n=2$, it follows
from \eqref{eq3.13} that estimate \eqref{eq3.3} is established. This
completes the proof of Lemma \ref{lem3.1}.
\end{proof}
%--------------------(end the proof of lem3.1)------------------------------------
\medskip

\textbf{Proof of Theorem \ref{thm1.1}:} By using standard method, we
only need to give the a priori estimates  to control
$\|u(t)\|_{H^{3}}+\|\nabla d(t)\|_{H^{3}}$ for any $0\leq t\leq
T_{*}$ in terms of $u_{0}$, $d_{0}$ and $\varepsilon_{0}$.  To this
end, we need to introduce the following commutator and product
estimates (see \cite{KP,BCD,PG}):
\begin{align}
   \label{eq3.14}
&\|\Lambda^{\alpha}(fg)-f\Lambda^{\alpha}g\|_{L^{p}}\leq C (\|\nabla
f\|_{L^{p_{1}}}\|\Lambda^{\alpha-1}g\|_{L^{q_{1}}}+\|\Lambda^{\alpha}f\|_{L^{p_{2}}}\|g\|_{L^{q_{2}}});\\
%--------------------(eq3.14)------------------------------------
   \label{eq3.15}
&\|\Lambda^{\alpha}(fg)\|_{L^{p}}\leq C
(\|f\|_{L^{p_{1}}}\|\Lambda^{\alpha}g\|_{L^{q_{2}}}+\|\Lambda^{\alpha}f\|_{L^{p_{2}}}\|g\|_{L^{q_{2}}})
\end{align}
%--------------------(eq3.15)------------------------------------
with $\alpha>0$, $1<p,p_{1},p_{2},q_{1},q_{2}<\infty$ and
$\frac{1}{p}=\frac{1}{p_{1}}+\frac{1}{q_{1}}=\frac{1}{p_{2}}+\frac{1}{q_{2}}$.
Here $\Lambda: =(-\Delta)^{\frac{1}{2}}$.

Applying $\Lambda^{3}$ on \eqref{eq1.1}, multiplying $\Lambda^{3}u$,
integrating with respect to $x$ over $\mathbb{R}^{3}$, and using
integration by parts, one obtains
\begin{align}\label{eq3.16}
\!\frac{1}{2}\frac{d}{dt}\|\Lambda^{3}u
(\cdot,t)\|_{\!L^{\!2}}^{2}\!+\!\|\Lambda^{4}u
(\cdot,t)\|_{\!L^{\!2}}^{2}\!=\!\!-\!\!\int_{\mathbb{R}^{3}}\!\!
\Lambda^{3}(u\cdot\nabla
u)\cdot\Lambda^{3}u\text{d}x\!-\!\!\!\int_{\mathbb{R}^{3}}\!\!\Lambda^{3}(\Delta
d\cdot\nabla
d)\cdot\Lambda^{3}u\text{d}x\!:=\!I_{\!6}\!\!+\!I_{\!7}.
\end{align}
%--------------------(eq3.16)------------------------------------
Noticing that the fact that $\operatorname{div} u=0$ implies
$\int_{\mathbb{R}^{3}}\Lambda^{3}\nabla(\frac{|\nabla
d|^{2}}{2})\cdot\Lambda^{3} u\text{d}x=0$, it follows that
\begin{align}\label{eq3.17}
I_{6}=&\int_{\mathbb{R}^{3}}[\Lambda^{3}(u\cdot\nabla
u)-u\cdot\nabla\Lambda^{3}u]\cdot\Lambda^{3}u\text{d}x\nonumber\\
\leq& C\|[\Lambda^{3}(u\cdot\nabla
u)-u\cdot\nabla\Lambda^{3}u]\|_{L^{\frac{3}{2}}}\|\Lambda^{3}u\|_{L^{3}}\nonumber\\
\leq &C\|\nabla u\|_{L^{3}}\|\Lambda^{3} u\|_{L^{3}}^{2}\leq
C\|\nabla
u\|_{L^{2}}^{\frac{7}{6}}\|\Lambda^{4}u\|_{L^{2}}^{\frac{11}{6}}\nonumber\\
\leq& \frac{1}{4}\|\Lambda^{4} u\|_{L^{2}}^{2}+C\|\nabla
u\|_{L^{2}}^{14}\nonumber\\
\leq& \frac{1}{4}\|\Lambda^{4} u\|_{L^{2}}^{2}+CC_{0}^{7},
\end{align}
%--------------------(eq3.17)------------------------------------
where $C_{0}$ is the bounded positive constant in \eqref{eq3.3}.
Here we have used the following Gagliardo--Nirenberg inequalities:
\begin{align*}
\|\nabla u\|_{L^{3}}\leq C\|\nabla
u\|_{L^{2}}^{\frac{5}{6}}\|\Lambda^{4} u\|_{L^{2}}^{\frac{1}{6}}
\text{ and }
 \|\Lambda^{3} u\|_{L^{3}}\leq C\|\nabla
u\|_{L^{2}}^{\frac{1}{6}}\|\Lambda^{4}u\|_{L^{2}}^{\frac{5}{6}}.
\end{align*}
For $I_{7}$, applying the H\"{o}lder's inequality and the Leibniz's
rule,  we have
\begin{align}\label{eq3.18}
I_{7}= &\int_{\mathbb{R}^{3}}\Lambda^{2}(\Delta d\cdot\nabla d)\cdot
\Lambda^{4}u\text{d}x\nonumber\\
\leq & \frac{1}{4}\|\Lambda^{4}
u\|_{L^{2}}^{2}+C\int_{\mathbb{R}^{3}}|\Lambda^{2}(\Delta
d\cdot\nabla d)|^{2}\text{d}x\nonumber\\
\leq &\frac{1}{4}\|\Lambda^{4}
u\|_{L^{2}}^{2}+C\int_{\mathbb{R}^{3}}(|\Lambda^{4}d|^{2}|\nabla
d|^{2}+|\Lambda^{2} d|^{2}|\Lambda^{3} d|^{2})\text{d}x\nonumber\\
\leq& \frac{1}{4}\|\Lambda^{4} u\|_{L^{2}}^{2}+C(\|\nabla
d\|_{L^{6}}^{2}\|\Lambda^{4} d\|_{L^{3}}^{2}+\|\Lambda^{2}d\|_{L^{4}}^{2}\|\Lambda^{3}d\|_{L^{4}}^{2})\nonumber\\
\leq& \frac{1}{4}\|\Lambda^{4} u\|_{L^{2}}^{2}+C(\|\Delta
d\|_{L^{2}}^{\frac{7}{3}}\|\Lambda^{5}d\|_{L^{2}}^{\frac{5}{3}}+
\|\Delta d\|_{L^{2}}^{\frac{19}{6}}\|\Lambda^{5}
d\|_{L^{2}}^{\frac{5}{6}})\nonumber\\
\leq & \frac{1}{4}\|\Lambda^{4}
u\|_{L^{2}}^{2}+\frac{1}{4}\|\Lambda^{5} d\|_{L^{2}}^{2} +C(
\|\Delta d\|_{L^{2}}^{14}+\|\Delta
d\|_{L^{2}}^{\frac{38}{7}})\nonumber\\
\leq & \frac{1}{4}\|\Lambda^{4}
u\|_{L^{2}}^{2}+\frac{1}{4}\|\Lambda^{5} d\|_{L^{2}}^{2}
+C(C_{0}^{7}+C_{0}^{\frac{19}{7}}).
\end{align}
%--------------------(eq3.18)------------------------------------
Here we have used the following Gagliardo--Nirenberg inequalities:
\begin{align*}
&\|\Lambda^{4} d\|_{L^{3}}\leq C\|\Delta
d\|_{L^{2}}^{\frac{1}{6}}\|\Lambda^{5} d\|_{L^{2}}^{\frac{5}{6}};\\
& \|\Lambda^{2} d\|_{L^{4}}\leq C\|\Delta
d\|_{L^{2}}^{\frac{3}{4}}\|\Lambda^{5}d\|_{L^{2}}^{\frac{1}{4}};\\
&\|\Lambda^{3}d\|_{L^{4}}\leq C\|\Delta
d\|_{L^{2}}^{\frac{5}{6}}\|\Lambda^{5}d\|_{L^{2}}^{\frac{1}{6}}.
\end{align*}
Inserting \eqref{eq3.17} and \eqref{eq3.18} into \eqref{eq3.16}, one
gets
\begin{align}\label{eq3.19}
\frac{d}{dt}\|\Lambda^{3} u\|_{L^{2}}^{2}+&\|\Lambda^{4}
u\|_{L^{2}}^{2}\leq
\frac{1}{2}\|\Lambda^{5}d\|_{L^{2}}^{2}+C(C_{0}^{7}+C_{0}^{\frac{19}{7}}).
\end{align}
%--------------------(eq3.19)------------------------------------

Taking $\Lambda^{4}$ on \eqref{eq1.2}, multiplying $\Lambda^{4} d$,
integrating with respect to $x$ over $\mathbb{R}^{3}$, and using
integration by parts, one obtains
\begin{align}\label{eq3.20}
\frac{1}{2}\frac{d}{dt}\|\Lambda^{4} d\|_{L^{2}}^{2}+\|\Lambda^{5}
d\|_{L^{2}}^{2}=-\int_{\mathbb{R}^{3}}\Lambda^{4}(u\cdot\nabla
d)\cdot\Lambda^{4}d\text{d}x+\int_{\mathbb{R}^{3}}\Lambda^{4}(|\nabla
d|^{2}d)\cdot\Lambda^{4}d\text{d}x:=I_{8}+I_{9}.
\end{align}
%--------------------(eq3.20)------------------------------------
Similar as estimate of $I_{6}$, we have
\begin{align}\label{eq3.21}
I_{8}=&-\int_{\mathbb{R}^{3}}[\Lambda^{4} (u\cdot\nabla
d)-u\cdot\nabla\Lambda^{4} d]\cdot \Lambda^{4}d\text{d}x\nonumber\\
\leq&  C\|\Lambda^{4} (u\cdot\nabla d)-u\cdot\nabla\Lambda^{4}
d\|_{L^{\frac{3}{2}}}\|\Lambda^{4} d\|_{L^{3}}\nonumber\\
\leq& C\|\nabla d\|_{L^{6}}\|\Lambda^{4} u\|_{L^{2}}\|\Lambda^{4}
d\|_{L^{3}}+C\|\nabla u\|_{L^{2}}\|\Lambda^{4} d\|_{L^{2}}\|\Lambda^{4} d\|_{L^{3}}\nonumber\\
\leq& \frac{1}{4}\|\Lambda^{4} u\|_{L^{2}}^{2} +C\|\Delta
d\|_{L^{2}}^{2}\|\Lambda^{4} d\|_{L^{3}}^{2}+\|\nabla u\|_{L^{2}}\|\Lambda^{5} d\|_{L^{2}}\|\Lambda^{4} d\|_{L^{3}} \nonumber\\
\leq& \frac{1}{4}\|\Lambda^{4} u\|_{L^{2}}^{2} +C\|\Delta
d\|_{L^{2}}^{2}\|\Delta d\|_{L^{2}}^{\frac{1}{3}}\|\Lambda^{5}
d\|_{L^{2}}^{\frac{5}{3}}
+\|\nabla u \|_{L^{2}}\|\Delta d\|_{L^{2}}^{\frac{1}{6}}\|\Lambda^{5} d\|_{L^{2}}^{\frac{11}{6}}\nonumber\\
\leq & \frac{1}{4}\|\Lambda^{4} u\|_{L^{2}}^{2}+\frac{1}{4}
\|\Lambda^{5} d\|_{L^{2}}^{2} +C (\|\Delta d\|_{L^{2}}^{14}+\|\nabla
u\|_{L^{2}}^{24}+\|\Delta d\|_{L^{2}}^{4})\nonumber\\
\leq & \frac{1}{4}\|\Lambda^{4} u\|_{L^{2}}^{2}+\frac{1}{4}
\|\Lambda^{5} d\|_{L^{2}}^{2} +C (C_{0}^{7}+C_{0}^{12}+C_{0}^{2}),
\end{align}
%--------------------(eq3.21)------------------------------------
where we have used the  Gagliardo--Nirenberg inequality:
\begin{align*}
\|\Lambda^{4} d\|_{L^{2}} \leq C\|\Delta
d\|_{L^{2}}^{\frac{1}{3}}\|\Lambda^{5}
d\|_{L^{2}}^{\frac{2}{3}}\text{ and }
 \|\Lambda^{4} d\|_{L^{3}} \leq
C\|\Delta d\|_{L^{2}}^{\frac{1}{6}}\|\Lambda^{5}
d\|_{L^{2}}^{\frac{5}{6}}.
\end{align*}
To estimate $I_{9}$, by using the Leibniz's rule, the fact $|d|=1$,
the H\"{o}lder's inequality               %, the interpolation inequalities
 and the
Young inequality, one obtains
\begin{align}\label{eq3.22}
I_{9}=&\int_{\mathbb{R}^{3}}\Lambda^{4}(|\nabla d|^{2}d)\cdot
\Lambda^{4}d\text{d}x=-\int_{\mathbb{R}^{3}}\Lambda^{3}(|\nabla
d|^{2}d)\cdot \Lambda^{5}d\text{d}x\nonumber\\
=&-\!\int_{\mathbb{R}^{3}}\!\big[\Lambda^{3}(|\nabla d|^{2})
d\cdot\Lambda^{5}d+3\Lambda^{2}(|\nabla d|^{2})\Lambda d
\cdot\Lambda^{5} d +3\Lambda (|\nabla d|^{2})\Lambda^{2} d \cdot
\Lambda^{5} d +|\nabla d|^{2}\Lambda^{3}d\cdot\Lambda^{5}
d\big]\text{d}x\nonumber\\
\leq & C\|\Lambda^{5} d\|_{L^{2}}(\|\nabla d\|_{L^{6}}\|\Lambda^{4}
d\|_{L^{3}}+\|\Lambda^{2}
d\|_{L^{4}}\|\Lambda^{3}d\|_{L^{4}}+\|\nabla
d\|_{L^{6}}^{2}\|\Lambda^{3} d\|_{L^{6}}+\|\nabla
d\|_{L^{6}}\|\Lambda^{2} d\|_{L^{6}}^{2})\nonumber\\
\leq & C\|\Lambda^{5} d\|_{L^{2}}(\|\Delta
d\|_{L^{2}}^{\frac{7}{6}}\|\Lambda^{5}d\|_{L^{2}}^{\frac{5}{6}} +\|
\Delta d\|_{L^{2}}^{\frac{7}{3}}\|\Lambda^{5}
d\|_{L^{2}}^{\frac{2}{3}})\nonumber\\
\leq &  \frac{1}{4}\|\Lambda^{5} d\|_{L^{2}}^{2}+C\|\Delta
d\|_{L^{2}}^{14}\nonumber\\
\leq &  \frac{1}{4}\|\Lambda^{5} d\|_{L^{2}}^{2}+C C_{0}^{7}.
\end{align}
%--------------------(eq3.22)------------------------------------
Here we have used the following Gagliardo--Nirenberg inequalities:
\begin{align*}
&\|\Lambda^{4} d\|_{L^{3}}\leq C\|\Delta
d\|_{L^{2}}^{\frac{1}{6}}\|\Lambda^{5}
d\|_{L^{2}}^{\frac{5}{6}};\quad\quad
 \|\Lambda^{2} d\|_{L^{4}}\leq
C\|\Delta
d\|_{L^{2}}^{\frac{3}{4}}\|\Lambda^{5}d\|_{L^{2}}^{\frac{1}{4}};\\
&\|\Lambda^{3}d\|_{L^{6}}\leq C\|\Delta
d\|_{L^{2}}^{\frac{1}{3}}\|\Lambda^{5}d\|_{L^{2}}^{\frac{2}{3}};\quad\quad
\|\Lambda^{2} d\|_{L^{6}}\leq C \|\Delta
d\|_{L^{2}}^{\frac{2}{3}}\|\Lambda^{5} d\|_{L^{2}}^{\frac{1}{3}}.
\end{align*}
Inserting \eqref{eq3.21} and \eqref{eq3.22} into \eqref{eq3.20}, one
gets
\begin{align}\label{eq3.23}
\frac{d}{dt}\|\Lambda^{4} d\|_{L^{2}}^{2}+\|\Lambda^{5}
d\|_{L^{2}}^{2}\leq& \frac{1}{2} \|\Lambda^{4} u\|_{L^{2}}^{2}++C
(C_{0}^{7}+C_{0}^{12}+C_{0}^{2}).
\end{align}
%--------------------(eq3.23)------------------------------------
Combining \eqref{eq3.19} and \eqref{eq3.23} together, and letting
$C_{0}>1$, one obtains
\begin{align*}
&\frac{d}{dt}(\|\Lambda^{3} u\|_{L^{2}}^{2}+\|\Lambda^{4}
d\|_{L^{2}}^{2})+\frac{1}{2}(\|\Lambda^{4}
u\|_{L^{2}}^{2}+\|\Lambda^{5} d\|_{L^{2}}^{2}) \leq C C_{0}^{14}.
\end{align*}
Hence integrating with respect to $t$ over $[0,T_{*}]$, we have
\begin{align*}
\sup_{0<t\leq T_{*}}(\|\Lambda^{3} u (t)\|_{L^{2}}^{2}+\|\Lambda^{4}
d (t)\|_{L^{2}}^{2})+\frac{1}{2}\int_{0}^{T_{*}}\left(\|\Lambda^{4}
u(\cdot,\tau)\|_{L^{2}}^{2}+\|\Lambda^{5}
d(\cdot,\tau)\|_{L^{2}}^{2}\right)\text{d}\tau\leq C<\infty,
\end{align*}
where $C$ only depends on the initial data $(u_{0},d_{0})$, $C_{0}$
and  $T_{*}$. Therefore, we get
\begin{align*}
&\|u\|_{L^{\infty}(0,T_{*};H^{3})}+\|u\|_{L^{2}(0,T_{*};H^{4})}\leq
C<\infty,\nonumber\\
&\|d\|_{L^{\infty}(0,T_{*};H^{4})}+\|d\|_{L^{2}(0,T_{*};H^{5})}\leq
C<\infty.
\end{align*}
This completes the proof of Theorem \ref{thm1.1}.  $\hfill\Box$
\medskip
\medskip

\textbf{Proof of Corollary \ref{cor1.2}:}  Similar as the proof of
Theorem \ref{thm1.1},  we only give the a priori estimates  to
control $\|u(t)\|_{H^{2}}+\|\nabla d(t)\|_{H^{3}}$ for any $0\leq
t\leq T_{*}$ in terms of $u_{0}$, $d_{0}$ and $\varepsilon_{0}$. To
this end, let us firstly recall the following useful
Gagliardo-Nirenberg inequalities in $\mathbb{R}^{2}$:
\begin{align}\label{eq3.24}
&\|\nabla u\|_{L^{3}}\leq C\|\nabla
u\|_{L^{2}}^{\frac{5}{6}}\|\Lambda^{3} u\|_{L^{2}}^{\frac{1}{6}};\quad\quad %\nonumber\\
\|\nabla u\|_{L^{6}}\leq C\|
u\|_{L^{2}}^{\frac{4}{9}}\|\Lambda^{3} u\|_{L^{2}}^{\frac{5}{9}};\nonumber\\
&\|\Lambda^{2} u\|_{L^{3}}\leq C\|\nabla
u\|_{L^{2}}^{\frac{1}{3}}\|\Lambda^{3} u\|_{L^{2}}^{\frac{2}{3}};\quad\quad %\nonumber\\
\|\nabla d\|_{L^{6}}\leq C\|\nabla
d\|_{L^{2}}^{\frac{1}{3}}\|\Delta d\|_{L^{2}}^{\frac{2}{3}};\\
&\|\Lambda^{2} d\|_{L^{4}}\leq C\|\Delta
d\|_{L^{2}}^{\frac{3}{4}}\|\Lambda^{4} d\|_{L^{2}}^{\frac{1}{4}};\quad\quad %\nonumber\\
\|\Lambda^{3} d\|_{L^{3}}\leq C\|\Delta
d\|_{L^{2}}^{\frac{1}{3}}\|\Lambda^{4}
d\|_{L^{2}}^{\frac{2}{3}}.\nonumber
\end{align}
%--------------------(eq3.24)------------------------------------
Now, applying $\Lambda^{2}$ on \eqref{eq1.1}, multiplying
$\Lambda^{2} u$ and integrating with respect to $x$ over
$\mathbb{R}^{3}$, and using \eqref{eq3.14}, the H\"{o}lder's
inequality, \eqref{eq3.24} and the Young inequality, one obtains
\begin{align}\label{eq3.25}
&\frac{1}{2}\frac{d}{dt}\|\Lambda^{2}u(\cdot,t)\|_{L^{2}}^{2}+\|\Lambda^{3}u(\cdot,t)\|_{L^{2}}^{2}
=-\int_{\mathbb{R}^{2}}\Lambda^{2}(u\cdot\nabla
u)\cdot\Lambda^{2}u\text{d}x-\int_{\mathbb{R}^{2}}\Lambda^{2}(\Delta
d\cdot \nabla d)\cdot\Lambda^{2}u\text{d}x\nonumber\\
=&-\int_{\mathbb{R}^{2}}[\Lambda^{2}(u\cdot\nabla u)-u\cdot\nabla
\Lambda^{2} u]\cdot\Lambda^{2}
u\text{d}x+\int_{\mathbb{R}^{2}}\Lambda(\Delta
d\cdot \nabla d)\cdot\Lambda^{3}u\text{d}x\nonumber\\
\leq & C\|[\Lambda^{2}(u\cdot\nabla u)-u\cdot\nabla \Lambda^{2}
u]\|_{L^{\frac{3}{2}}}\|\Lambda^{2} u\|_{L^{3}}+C(\|\nabla
d\|_{L^{6}}\|\Lambda^{3}d\|_{L^{3}}\|\Lambda^{3} u\|_{L^{2}}
+\|\Lambda^{3}u\|_{L^{2}}\|\Lambda^{2} d\|_{L^{4}}^{2})\nonumber\\
\leq& C (\|\Lambda^{2} u\|_{L^{3}}^{2}\|\nabla u\|_{L^{3}}+\|\nabla
d\|_{L^{6}}\|\Lambda^{3}d\|_{L^{3}}\|\Lambda^{3} u\|_{L^{2}}
+\|\Lambda^{3}u\|_{L^{2}}\|\Lambda^{2} d\|_{L^{4}}^{2})\nonumber\\
\leq &C (\|\nabla
u\|_{\!L^{\!2}}^{\frac{3}{2}}\|\Lambda^{3}u\|_{\!L^{\!2}}^{\frac{3}{2}}\!+\|\Lambda^{3}
u\|_{\!L^{\!2}}\|\nabla d\|_{\!L^{\!2}}^{\frac{1}{3}}\|\Delta
d\|_{\!L^{2}}\|\Lambda^{4}
d\|_{\!L^{\!2}}^{\frac{2}{3}}\!+\|\Lambda^{3}
u\|_{\!L^{\!2}}\|\Delta d\|_{\!L^{\!2}}^{\frac{3}{2}}\|\Lambda^{4}
d\|_{\!L^{\!2}}^{\frac{1}{2}})\nonumber\\
\leq & \frac{1}{4}\|\Lambda^{3}
u\|_{L^{2}}^{2}+\frac{1}{4}\|\Lambda^{4} d\|_{L^{2}}^{2}+C(\|\nabla
u\|_{L^{2}}^{6}+ \|\nabla d\|_{L^{2}}^{2}\|\Delta
d\|_{L^{2}}^{6}+\|\Delta d\|_{L^{2}}^{6})\nonumber\\
\leq & \frac{1}{4}\|\Lambda^{3}
u\|_{L^{2}}^{2}+\frac{1}{4}\|\Lambda^{4}
d\|_{L^{2}}^{2}+C(1+\|\nabla u\|_{L^{2}}^{6}+ \|\Delta
d\|_{L^{2}}^{6})\nonumber\\
\leq & \frac{1}{4}\|\Lambda^{3}
u\|_{L^{2}}^{2}+\frac{1}{4}\|\Lambda^{4}
d\|_{L^{2}}^{2}+C(1+C_{0}^{3}),
\end{align}
%--------------------(eq3.25)------------------------------------
where we have used the energy equality \eqref{eq3.4}, and $C_{0}$ is
the positive constant defined in Lemma \ref{eq3.1}.

Taking $\Lambda^{3}$ on \eqref{eq1.2}, multiplying $\Lambda^{3} d$,
integrating with respect to $x$ over $\mathbb{R}^{2}$, and using
\eqref{eq3.14}, the H\"{o}lder's inequality, \eqref{eq3.24} and the
Young inequality, one obtains
\begin{align}\label{eq3.26}
\frac{1}{2}&\frac{d}{dt}\|\Lambda^{3}d(\cdot,t)\|_{L^{2}}^{2}+\|\Lambda^{4}
d
(\cdot,t)\|_{L^{2}}^{2}=-\int_{\mathbb{R}^{2}}\Lambda^{3}(u\cdot\nabla
d)\cdot
\Lambda^{3}d\text{d}x+\int_{\mathbb{R}^{2}}\Lambda^{3}(|\nabla
d|^{2}d)\cdot\Lambda^{3}d\text{d}x\nonumber\\
=& -\int_{\mathbb{R}^{2}}[\Lambda^{3}(u\cdot\nabla
d)-u\cdot\nabla\Lambda^{3}d]\cdot
\Lambda^{3}d\text{d}x-\int_{\mathbb{R}^{2}}\Lambda^{2}(|\nabla
d|^{2}d)\cdot\Lambda^{4}d\text{d}x\nonumber\\
=& -\!\!\int_{\mathbb{R}^{2}}\![\Lambda^{3}(u\cdot\nabla
d)-u\cdot\nabla\Lambda^{3}d]\cdot
\Lambda^{3}d\text{d}x\!-\!\!\int_{\mathbb{R}^{2}}\![\Lambda^{2}(|\nabla
d|^{2})d+\!2\Lambda(|\nabla d|^{2})\Lambda d
+\!|\nabla d|^{2}\Lambda^{2} d]\cdot\Lambda^{4}d\text{d}x\nonumber\\
\leq& C\|[\Lambda^{3}(u\cdot\nabla
d)-u\cdot\nabla\Lambda^{3}d]\|_{L^{\frac{3}{2}}}\|\Lambda^{3}d\|_{L^{3}}
+C(\|\Lambda^{3} d\|_{L^{3}}\|\nabla d\|_{L^{6}}+\|\Lambda^{2}
d\|_{L^{4}}^{2}\nonumber\\
&+\|\Lambda d\|_{L^{6}}^{2}\|\Lambda^{2} d\|_{L^{6}}+\|\nabla
d\|_{L^{6}}\|\Lambda^{3}d\|_{L^{3}})\|\Lambda^{4}
d\|_{L^{2}}\nonumber\\
\leq &C (\|\Lambda^{3} d\|_{L^{3}}\|\nabla d\|_{L^{6}}\|\Lambda^{3}
u\|_{L^{2}}+\|\Lambda^{3} d\|_{L^{3}}^{2}\|\nabla u\|_{L^{6}}+
\|\Lambda^{3} d\|_{L^{3}}\|\nabla d\|_{L^{6}}\|\Lambda^{4}
d\|_{L^{2}}\nonumber\\
&+\|\Lambda^{2} d\|_{L^{4}}^{2}\|\Lambda^{4} d\|_{L^{2}}+\|\Lambda
d\|_{L^{6}}^{2}\|\Lambda^{2} d\|_{L^{6}}\|\Lambda^{4}
d\|_{L^{2}})\nonumber\\
\leq & \frac{1}{8}\|\Lambda^{3}
u\|_{L^{2}}^{2}+\frac{1}{8}\|\Lambda^{4} d\|_{L^{2}}^{2} +C(\|\nabla
d\|_{L^{6}}^{2}\|\Lambda^{3} d\|_{L^{3}}^{2}+\|\nabla
u\|_{L^{6}}\|\Lambda^{3} d\|_{L^{3}}^{2}+\|\Lambda^{2}
d\|_{L^{4}}^{4}+\|\Lambda d\|_{L^{6}}^{4}\|\Lambda^{2}
d\|_{L^{6}}^{2})\nonumber\\
\leq & \frac{1}{8}\|\Lambda^{3}
u\|_{L^{2}}^{2}+\frac{1}{8}\|\Lambda^{4} d\|_{L^{2}}^{2} +C(\|\nabla
d\|_{L^{2}}^{\frac{2}{3}}\|\Delta d\|_{L^{2}}^{2}\|\Lambda^{4}
d\|_{L^{2}}^{\frac{4}{3}} +\|u\|_{L^{2}}^{\frac{4}{9}} \|\Lambda^{3}
u\|_{L^{2}}^{\frac{5}{9}} \|\Delta
d\|_{L^{2}}^{\frac{2}{3}}\|\Lambda^{4}
d\|_{L^{2}}^{\frac{4}{3}}\nonumber\\
& +\|\Delta d\|_{L^{2}}^{3}\|\Lambda^{4} d\|_{L^{2}}+\|\nabla
d\|_{L^{2}}^{\frac{4}{3}}\|\Delta d\|_{L^{2}}^{4}\|\Lambda^{4}
d\|_{L^{2}}^{\frac{2}{3}})\nonumber\\
\leq & \frac{1}{4}\|\Lambda^{3}
u\|_{L^{2}}^{2}+\frac{1}{4}\|\Lambda^{4} d\|_{L^{2}}^{2} +
C(\|\nabla d\|_{L^{2}}^{2}\|\Delta
d\|_{L^{2}}^{6}+\|u\|_{L^{2}}^{8}\|\Delta d\|_{L^{2}}^{12}+\|\Delta
d\|_{L^{2}}^{6}+\|\nabla d\|_{L^{2}}^{2}\|\Delta d\|_{L^{2}}^{6}
)\nonumber\\
\leq & \frac{1}{4}\|\Lambda^{3}
u\|_{L^{2}}^{2}+\frac{1}{4}\|\Lambda^{4} d\|_{L^{2}}^{2} +
C(\|\Delta d\|_{L^{2}}^{6}+\|\Delta
d\|_{L^{2}}^{12})\nonumber\\
\leq& \frac{1}{4}\|\Lambda^{3}
u\|_{L^{2}}^{2}+\frac{1}{4}\|\Lambda^{4} d\|_{L^{2}}^{2} +
C(1+\|\Delta d\|_{L^{2}}^{12})\nonumber\\
\leq& \frac{1}{4}\|\Lambda^{3}
u\|_{L^{2}}^{2}+\frac{1}{4}\|\Lambda^{4} d\|_{L^{2}}^{2} +
C(1+C_{0}^{6}),
\end{align}
%--------------------(eq3.26)------------------------------------
where we have used the energy equality \eqref{eq3.4}, and $C_{0}$
defined in Lemma \ref{eq3.1}. Combining \eqref{eq3.25} and
\eqref{eq3.26} together,  we obtain
\begin{align}\label{eq3.27}
&\frac{d}{dt}(\|\Lambda^{2}u(\cdot,t)\|_{L^{2}}^{2}+\|\Lambda^{3}d(\cdot,t)\|_{L^{2}}^{2})
+(\|\Lambda^{3}u\|_{L^{2}}^{2} +\|\Lambda^{4} d \|_{L^{2}}^{2})\leq
C(1+C_{0}^{6}).
\end{align}
%--------------------(eq3.27)------------------------------------
Hence integrating with respect to $t$ over $[0,T_{*}]$, we have
\begin{align*}
\sup_{0<t\leq T_{*}}(\|\Lambda^{2} u (t)\|_{L^{2}}^{2}+\|\Lambda^{3}
d (t)\|_{L^{2}}^{2})+\int_{0}^{T_{*}}\left(\|\Lambda^{3}
u(\cdot,\tau)\|_{L^{2}}^{2}+\|\Lambda^{4}
d(\cdot,\tau)\|_{L^{2}}^{2}\right)\text{d}\tau\leq C<\infty,
\end{align*}
where $C$ only depends on the initial data $(u_{0},d_{0})$, $C_{0}$
and $T_{*}$. Therefore, we get
\begin{align*}
&\|u\|_{L^{\infty}(0,T_{*};H^{2})}+\|u\|_{L^{2}(0,T_{*};H^{3})}\leq
C<\infty,\nonumber\\
&\|d\|_{L^{\infty}(0,T_{*};H^{3})}+\|d\|_{L^{2}(0,T_{*};H^{4})}\leq
C<\infty.
\end{align*}
This completes the proof of Corollary \ref{cor1.2}.  $\hfill\Box$
\\
\\

%--------------------(proof of thm1.4)------------------------------------
%\\
%\\
%\textbf{Acknowledgments}


\begin{thebibliography}{99}


\bibitem{BCD} H. Bahouri, J. Chemin and R. Danchin, \textit{Fourier
 Analysis and Nonlinear Partial Differential Equations}, Springer
 Heidelberg Dordrecht London New York, 2011.

\bibitem{HBV} H. Beir\~{a}o da Veiga, A new regularity class for the
Navier--Stokes equations in $\mathbb{R}^{n}$, Chinese Ann. Math.
Ser., B 16 (1995), 407--412.

\bibitem{BKM} J. Beale, T. Kato and A. Majda, Remarks on breakdown of smooth solutions for
the 3D Euler equations. Commun. Math. Phys., 94 (1984), 61--66.

\bibitem{CPG} D. Charmorro and P. G. Lemari\'{e}-Rieusset, Real
Interpolation method, Lorentz spaces and refined Sobolev
inequalities, arXiv:1211.3320v1 [math.AP] 14 Nov 2012.

\bibitem{CDY} K. Chang, W. Ding and R. Ye, Finite-time blow-up of
the heat flow of harmonic maps from surfaces, J. Differ. Geom.,
36(2) (1992), 507--515.

\bibitem{CS} Y. Chen and M. Struwe, Existence and partial regularity results for the heat
 flow of harmonic maps, Math. Z., 201 (1989), 83--103.


\bibitem{ER} J. L. Ericksen, Hydrostatic theory of liquid crystal, Arch. Rational Mech. Anal.,
 9 (1962), 371--378.

\bibitem{ESS} L. Escauriaza, G. Seregin and V. \u{S}ver\'{a}k,
Backward uniqueness for parabolic equaitons, Arch. Rational Mech.
Anal., 169 (2003), 147--157.


\bibitem{YG} Y. Giga, Solutions for semilinear parabolic equations in $L^{p}$ and regularity of weak solutions of the
Navier--Stokes system, J. Differ. Equ., 61 (1986), 186--212.

\bibitem{GG11} Z. Guo and S. Gala, Remarks on logarithmical
regularity criteria for the Navier--Stokes equations, J, Math.
Phys., 52 (2011), 063503.



%\bibitem{HKL} R. Hardt, D. Kinderleher and F. H. Lin, Existence and
%partial regularity of static liquid crystal configurations, Commun.
%Math. Phys., 105 (1986), 547--570.

\bibitem{JHW} J. Hineman and C. Wang, Well--posedness of nematic
liquid crystal flow in $L^{3}_{loc}(\mathbb{R}^{3})$,
arXiv:1208.5965v1 [math.AP] 29 Aug. 2012.

\bibitem{HMC} M. Hong, Global existence of solutions of the
simplified Ericksen--Leslie system in dimension two, Cal. Var., 40
(2011), 15--36.


\bibitem{HW} X. Hu and D. Wang, Global Solution to the
Three-Dimensional Incompressible Flow of Liquid Crystals, Commun.
Math. Phys., 296 (2010), 861--880.

\bibitem{HW1} T. Huang and C. Wang, Blow up Criterion for Nematic
Liquid Crystal Flows, Comm. Partial Differ. Equ., 37 (2012),
875--884.

\bibitem{KP} T. Kato and G. Ponce, Commutator estimates and the
Euler and Navier--Stokes equations, Commu. Pure Appl. Math., 41
(1988), 891--907.

\bibitem{KOT} H. Kozono, T. Ogawa and Y. Taniuchi, The critical
Sobolev inequalities in Besov spaces and regularity criterion to
some semi-linear evolution equations, Math. Z., 242 (2002),
251--278.

\bibitem{KT} H. Kozono and Y. Taniuchi, Bilinear estimates in BMO
and the Navier--Stokes equations, Math. Z., 235 (2000), 173--194.

\bibitem{LE} F. Leslie, Theory of flow phenomenum in liquid crystals. In: The Theory of
Liquid Crystals, London-New York: Academic Press, 4 (1979), 1--81.

\bibitem{PG}  P.G. Lemari\'{e}-Rieusset, \textit{Recent Developments in the
 Navier--Stokes Problem}, Chapman and Hall/CRC, 2002.

\bibitem{XLW} X. Li and D. Wang, Global solution to the
incompressible flow of liquid crystal, J. Differ. Equ., 252 (2012),
745--767.

\bibitem{L} F. Lin, Nonlinear theory of defects in nematic liquid crystals; phase transition and flow phenomena,
 Comm. Pure. Appl. Math., 42 (1989), 789--814.

 \bibitem{LLW} F. Lin, J. Lin and C. Wang, Liquid Crystal flow in two dimensions, Arch.
 Rational Mech. Anal., 197 (2010), 297--336

\bibitem{LL1} F. Lin and C. Liu, Nonparabolic dissipative systems modeling the flow of liquid crystals,
Comm. Pure. Appl. Math.,  48  (1995), 501--537.

\bibitem{LL2} F. Lin and C. Liu, Partial regularities of the nonlinear dissipative systems modeling
the flow of liquid crystals, Disc. Contin. Dyn. Syst., A 2 (1996),
1--23.

\bibitem{LW} F. Lin and C. Wang, On the uniqueness of heat flow of harmonic maps and hydrodynamic
flow of nematic liquid crystals, Chinese Annal. Math. Ser., B 31
(2010), 921--938.

\bibitem{LD} J. Lin and S. Ding, On the well-posedness for the heat flow of harmonic maps and hydrodynamic
flow of nematic liquid crystals in critical spaces, Math. Meth.
Appl. Sciences, DOI: 10.1002/mma.1548.

\bibitem{LNW} C. Liu and N. J. Wakington, Approximation of Liquid
Crystal Flows, SIAM J. Numer. Anal., 37 (2000), 725--741.


%\bibitem{HM} H. Miura, Remark on uniqueness of mild solutions to the Navier--Stokes equations, J. Funct.
%Anal., 218 (2005), 110--129.

\bibitem{JS} J. Serrin, On the regularity of weak solutions of the
Navier--Stokes equations, Arch. Rational Mech. Anal., 9 (1962),
187--195.

%\bibitem{S} E. M. Stein, \textit{Singular integrals and differentiability properties of
%functions}, Princeton, NJ: Priceton University Press, 1971.

\bibitem{SL} H. Sun and C. Liu, On energetic variational approaches in modeling the nematic liquid crystal
flows, Disc. Contin. Dyn. Syst., A 23 (2009), 455--475.

%\bibitem{T} R. Temam, \textit{Navier--Stokes Equations}, North
%Holland, Amsterdam, 1977.

%\bibitem{TH} H. Triebel, \textit{Theory of Function Spaces.}
%Monograph in Mathematics, Vol. 78. Basel: Birkhauser Verglag, 1983.

\bibitem{W08} C. Wang, Heat flow of harmonic maps whose gradients
belong to $L^{n}_{x}L^{\infty}_{t}$, Arch. Rational Mech. Anal., 188
(2008), 309--349.


\bibitem{W} C. Wang, Well-posedness for the heat flow of harmonic
maps and the liquid crystal flow with rough initial data, Arch.
Rational Mech. Anal., 200 (2011), 1--19.

\bibitem{WD} H. Wen and S. Ding, Solutions of incompressible
hydrodynamic flow of liquid crystals, Nonlinear Anal. Real Word
Appl., 12 (2011), 1510--1531.




\end{thebibliography}
\end{document}